\documentclass{article}
\usepackage{amsmath,amsthm,amsopn,amstext,amscd,amsfonts,amssymb}
\usepackage{natbib}

\def \u {\mathop{\rm \mathcal{U}}\nolimits}
\def \tr {\mathop{\rm tr}\nolimits}
\def \re {\mathop{\rm Re}\nolimits}
\def \im {\mathop{\rm Im}\nolimits}
\def \E {\mathop{\rm E}\nolimits}

\def \Vol {\mathop{\rm Vol}\nolimits}
\def \Vec {\mathop{\rm vec}\nolimits}

\def \Cov {\mathop{\rm Cov}\nolimits}
\def \etr {\mathop{\rm etr}\nolimits}
\def \diag {\mathop{\rm diag}\nolimits}
\def \build#1#2#3{\mathrel{\mathop{#1}\limits^{#2}_{#3}}}

\renewenvironment{abstract}
                 {\vspace{6pt}
                  \begin{center}
                  \begin{minipage}{5in}
                  \centerline{\textbf{Abstract}}
                  \noindent\ignorespaces
                 }
                 {\end{minipage}\end{center}}
\newtheorem{theorem}{\textbf{Theorem}}[section]
\newtheorem{corollary}{\textbf{Corollary}}[section]
\newtheorem{lemma}{\textbf{Lemma}}[section]

\theoremstyle{definition}
\newtheorem{definition}{\textbf{Definition}}[section]

\newtheorem{remark}{\textbf{Remark}}[section]

\title{\Large \textbf{A generalised Kotz type distribution and Riesz distribution}}
\author{
  \textbf{Jos\'e A. D\'{\i}az-Garc\'{\i}a} \thanks{Corresponding author\newline
    }\\
  {\normalsize Department of Statistics and Computation} \\
  {\normalsize 25350 Buenavista, Saltillo, Coahuila, Mexico} \\
  {\normalsize E-mail: jadiaz@uaaan.mx} \\
}
\date{}
\begin{document}
\maketitle

\begin{abstract}
This article derives the distribution of random matrix $\mathbf{X}$ associated with the
transformation $\mathbf{Y} = \mathbf{X}^{*}\mathbf{X}$, such that $\mathbf{Y}$ has a Riesz
distribution for real normed division algebras. Two versions of this distributions are proposed
and some of their properties are studied.
\end{abstract}

\section{Introduction}\label{sec1}

Let $\mathbf{X}_{1}, \dots,\mathbf{X}_{n}$ be real independent $\mathcal{N}_{m}(\boldsymbol{\mu},
\boldsymbol{\Sigma})$, that is, $\mathbf{X}_{i}$ has a \emph{m-dimensional normal distribution} with
expected value $\boldsymbol{\mu} \in \Re^{m}$ and $m \times m$ positive definite covariance matrix
$\boldsymbol{\Sigma} > \mathbf{0}$. Let $\mathbf{X}$ be the $n \times m$ $(n \geq m)$, random matrix
$$
  \mathbf{X} =
  \left[
     \begin{array}{c}
       \mathbf{X}'_{1} \\
       \vdots \\
       \mathbf{X}'_{n}
     \end{array}
  \right]
$$
and observe that
$$
  \E(\mathbf{X}) =
  \left[
     \begin{array}{c}
       \boldsymbol{\mu}' \\
       \vdots \\
       \boldsymbol{\mu}'
     \end{array}
  \right] = \mathbf{1}\boldsymbol{\mu}', \quad \mbox{where } \mathbf{1} = (1, \dots,1)'\in \Re^{n}
$$
and $\Cov(\Vec \mathbf{X}') = \mathbf{I}_{n} \otimes \boldsymbol{\Sigma}$, that is, $\mathbf{X}$ has a $n
\times m$ \emph{matrix variate normal distribution}, denoted this fact as $\mathbf{X} \sim \mathcal{N}_{n
\times m}(\mathbf{1}\boldsymbol{\mu}', \mathbf{I}_{n}, \boldsymbol{\Sigma})$, see \citep[pp.
79-80]{m:82}, among others.

Now, if $\mathbf{Y} = \mathbf{X}'\mathbf{X}$, where the $n \times m$ random matrix $\mathbf{X}$ is
$\mathcal{N}_{n \times m}(\mathbf{0}, \mathbf{I}_{n}, \boldsymbol{\Sigma})$ then $\mathbf{Y}$ is said to
have the Wishart distribution with $n$ degrees of freedom and scale parameter $\boldsymbol{\Sigma}$.

Based on the theory of Jordan algebras, a family of distributions on symmetric cones, termed the
\emph{Riesz distributions},were first introduced by \citep{hl:01} under the name of Riesz natural
exponential family (Riesz NEF); they were based on a special case of the so-called \emph{Riesz measure}
from \citep[p.137]{fk:94}, going back to \citep{r:49}. This Riesz distribution \emph{generalises the
matrix multivariate gamma and Wishart distributions}, containing them as particular cases. Recently,
\citep{dg:15a} proposes two versions of the Riesz distribution and a diverse range of their properties
are studied for real normed division algebras.

The main purpose of this article is to introduce the distribution of the random matrix $\mathbf{X}$ such
that $\mathbf{Y} = \mathbf{X}^{*}\mathbf{X}$ have one of the two versions of the Riesz distributions and
explores some of their basic properties for real normed division algebras. These distributions shall be
termed  \emph{Kotz-Riesz distributions}. Moreover, as can be seen easily, the Kotz-Riesz distribution
belong to the \emph{matrix multivariate elliptical-spherical distributions}, see \citep[pp.
102-013]{fz:90}.

This article studies two versions of Kotz-Riesz distributions for real normed division algebras. Section
\ref{sec2} reviews some definitions and notation on real normed division algebras. And also, introduces
other mathematical tools as two definitions for the generalised gamma function on symmetric cones.
Several integration results for real normed division algebras are found in Section \ref{sec3}. Section
\ref{sec4} introduces Kotz-Riesz distributions for real normed division algebras, and also studies the
relationship between the Riesz distributions and the Kotz-Riesz distributions.

\section{Preliminary results}\label{sec2}

A detailed discussion of real normed division algebras may be found in \citep{b:02} and \citep{E:90}. For
your convenience, we shall introduce some notation, although in general we adhere to standard notation
forms.

For our purposes: Let $\mathbb{F}$ be a field. An \emph{algebra} $\mathfrak{A}$ over $\mathbb{F}$ is a
pair $(\mathfrak{A};m)$, where $\mathfrak{A}$ is a \emph{finite-dimensional vector space} over
$\mathbb{F}$ and \emph{multiplication} $m : \mathfrak{A} \times \mathfrak{A} \rightarrow A$ is an
$\mathbb{F}$-bilinear map; that is, for all $\lambda \in \mathbb{F},$ $x, y, z \in \mathfrak{A}$,
\begin{eqnarray*}
  m(x, \lambda y + z) &=& \lambda m(x; y) + m(x; z) \\
  m(\lambda x + y; z) &=& \lambda m(x; z) + m(y; z).
\end{eqnarray*}
Two algebras $(\mathfrak{A};m)$ and $(\mathfrak{E}; n)$ over $\mathbb{F}$ are said to be
\emph{isomorphic} if there is an invertible map $\phi: \mathfrak{A} \rightarrow \mathfrak{E}$ such that
for all $x, y \in \mathfrak{A}$,
$$
  \phi(m(x, y)) = n(\phi(x), \phi(y)).
$$
By simplicity, we write $m(x; y) = xy$ for all $x, y \in \mathfrak{A}$.

Let $\mathfrak{A}$ be an algebra over $\mathbb{F}$. Then $\mathfrak{A}$ is said to be
\begin{enumerate}
  \item \emph{alternative} if $x(xy) = (xx)y$ and $x(yy) = (xy)y$ for all $x, y \in \mathfrak{A}$,
  \item \emph{associative} if $x(yz) = (xy)z$ for all $x, y, z \in \mathfrak{A}$,
  \item \emph{commutative} if $xy = yx$ for all $x, y \in \mathfrak{A}$, and
  \item \emph{unital} if there is a $1 \in \mathfrak{A}$ such that $x1 = x = 1x$ for all $x \in \mathfrak{A}$.
\end{enumerate}
If $\mathfrak{A}$ is unital, then the identity 1 is uniquely determined.

An algebra $\mathfrak{A}$ over $\mathbb{F}$ is said to be a \emph{division algebra} if $\mathfrak{A}$ is
nonzero and $xy = 0_{\mathfrak{A}} \Rightarrow x = 0_{\mathfrak{A}}$ or $y = 0_{\mathfrak{A}}$ for all
$x, y \in \mathfrak{A}$.

The term ``division algebra", comes from the following proposition, which shows that, in such an algebra,
left and right division can be unambiguously performed.

Let $\mathfrak{A}$ be an algebra over $\mathbb{F}$. Then $\mathfrak{A}$ is a division algebra if, and
only if, $\mathfrak{A}$ is nonzero and for all $a, b \in \mathfrak{A}$, with $b \neq 0_{\mathfrak{A}}$,
the equations $bx = a$ and $yb = a$ have unique solutions $x, y \in \mathfrak{A}$.

In the sequel we assume $\mathbb{F} = \Re$ and consider classes of division algebras over $\Re$ or
``\emph{real division algebras}" for short.

We introduce the algebras of \emph{real numbers} $\Re$, \emph{complex numbers} $\mathfrak{C}$,
\emph{quaternions} $\mathfrak{H}$ and \emph{octonions} $\mathfrak{O}$. Then, if $\mathfrak{A}$ is an
alternative real division algebra, then $\mathfrak{A}$ is isomorphic to $\Re$, $\mathfrak{C}$,
$\mathfrak{H}$ or $\mathfrak{O}$.

Let $\mathfrak{A}$ be a real division algebra with identity $1$. Then $\mathfrak{A}$ is said to be
\emph{normed} if there is an inner product $(\cdot, \cdot)$ on $\mathfrak{A}$ such that
$$
  (xy, xy) = (x, x)(y, y) \qquad \mbox{for all } x, y \in \mathfrak{A}.
$$
If $\mathfrak{A}$ is a \emph{real normed division algebra}, then $\mathfrak{A}$ is isomorphic $\Re$,
$\mathfrak{C}$, $\mathfrak{H}$ or $\mathfrak{O}$.

There are exactly four normed division algebras: real numbers ($\Re$), complex numbers ($\mathfrak{C}$),
quaternions ($\mathfrak{H}$) and octonions ($\mathfrak{O}$), see \citep{b:02}. We take into account that,
$\Re$, $\mathfrak{C}$, $\mathfrak{H}$ and $\mathfrak{O}$ are the only normed division algebras;
furthermore, they are the only alternative division algebras.

Let $\mathfrak{A}$ be a division algebra over the real numbers. Then $\mathfrak{A}$ has dimension
$\beta$, either 1, 2, 4 or 8. Other branches of mathematics used the parameters $\alpha = 2/\beta$ and $t
= \beta/4$, see \citep{er:05} and \citep{k:84}, respectively.

Finally, observe that

\begin{tabular}{c}
  $\Re$ is a real commutative associative normed division algebras, \\
  $\mathfrak{C}$ is a commutative associative normed division algebras,\\
  $\mathfrak{H}$ is an associative normed division algebras, \\
  $\mathfrak{O}$ is an alternative normed division algebras. \\
\end{tabular}

Let $\mathfrak{L}^{\beta}_{n,m}$ be the set of all $n \times m$ matrices of rank $m \leq n$ over
$\mathfrak{A}$ with $m$ distinct positive singular values, where $\mathfrak{A}$ denotes a \emph{real
finite-dimensional normed division algebra}. Let $\mathfrak{A}^{n \times m}$ be the set of all $n \times
m$ matrices over $\mathfrak{A}$. The dimension of $\mathfrak{A}^{m \times n}$ over $\Re$ is $\beta mn$.
Let $\mathbf{A} \in \mathfrak{A}^{n \times m}$, then $\mathbf{A}^{*} = \bar{\mathbf{A}}^{T}$ denotes the
usual conjugate transpose.

Table \ref{table1} sets out the equivalence between the same concepts in the four normed division
algebras.

\begin{table}[th]
  \centering
  \caption{\scriptsize Notation}\label{table1}
  \begin{footnotesize}
  \begin{tabular}{cccc|c}
    \hline
    Real & Complex & Quaternion & Octonion & \begin{tabular}{c}
                                               Generic \\
                                               notation \\
                                             \end{tabular}\\
    \hline
    Semi-orthogonal & Semi-unitary & Semi-symplectic & \begin{tabular}{c}
                                                         Semi-exceptional \\
                                                         type \\
                                                       \end{tabular}
      & $\mathcal{V}_{m,n}^{\beta}$ \\
    Orthogonal & Unitary & Symplectic & \begin{tabular}{c}
                                                         Exceptional \\
                                                         type \\
                                                       \end{tabular} & $\mathfrak{U}^{\beta}(m)$ \\
    Symmetric & Hermitian & \begin{tabular}{c}
                              Quaternion \\
                              hermitian \\
                            \end{tabular}
     & \begin{tabular}{c}
                              Octonion \\
                              hermitian \\
                            \end{tabular} & $\mathfrak{S}_{m}^{\beta}$ \\
    \hline
  \end{tabular}
  \end{footnotesize}
\end{table}

It is denoted by ${\mathfrak S}_{m}^{\beta}$ the real vector space of all $\mathbf{S} \in \mathfrak{A}^{m
\times m}$ such that $\mathbf{S} = \mathbf{S}^{*}$. In addition, let $\mathfrak{P}_{m}^{\beta}$ be the
\emph{cone of positive definite matrices} $\mathbf{S} \in \mathfrak{A}^{m \times m}$. Thus,
$\mathfrak{P}_{m}^{\beta}$ consist of all matrices $\mathbf{S} = \mathbf{X}^{*}\mathbf{X}$, with
$\mathbf{X} \in \mathfrak{L}^{\beta}_{m,n}$; then $\mathfrak{P}_{m}^{\beta}$ is an open subset of
${\mathfrak S}_{m}^{\beta}$.

Let $\mathfrak{D}_{m}^{\beta}$ be the \emph{diagonal subgroup} of $\mathcal{L}_{m,m}^{\beta}$ consisting
of all $\mathbf{D} \in \mathfrak{A}^{m \times m}$, $\mathbf{D} = \diag(d_{1}, \dots,d_{m})$. Let
$\mathfrak{T}_{U}^{\beta}(m)$ be the subgroup of all \emph{upper triangular} matrices $\mathbf{T} \in
\mathfrak{A}^{m \times m}$ such that $t_{ij} = 0$ for $1 \leq j < i\leq m$.

The set of matrices $\mathbf{H}_{1} \in \mathfrak{F}^{n \times m}$ such that
$\mathbf{H}_{1}^{*}\mathbf{H}_{1} = \mathbf{I}_{m}$ is a manifold denoted ${\mathcal V}_{m,n}^{\beta}$,
is termed the \emph{Stiefel manifold} ($\mathbf{H}_{1}$ is also known as \emph{semi-orthogonal}).

For any matrix $\mathbf{X} \in \mathfrak{A}^{n \times m}$, $d\mathbf{X}$ denotes the\emph{ matrix of
differentials} $(dx_{ij})$. Finally, we define the \emph{measure} or volume element $(d\mathbf{X})$ when
$\mathbf{X} \in \mathfrak{A}^{n \times m}, \mathfrak{S}_{m}^{\beta}$, $\mathfrak{D}_{m}^{\beta}$ or
$\mathcal{V}_{m,n}^{\beta}$, see \citep{dggj:11}.

If $\mathbf{X} \in \mathfrak{A}^{n \times m}$ then $(d\mathbf{X})$ (the Lebesgue measure in
$\mathfrak{A}^{n \times m}$) denotes the exterior product of the $\beta mn$ functionally independent
variables
$$
  (d\mathbf{X}) = \bigwedge_{i = 1}^{m}\bigwedge_{j = 1}^{n}dx_{ij} \quad \mbox{ where }
    \quad dx_{ij} = \bigwedge_{k = 1}^{\beta}dx_{ij}^{(k)}.
$$

If $\mathbf{S} \in \mathfrak{S}_{m}^{\beta}$ (or $\mathbf{S} \in \mathfrak{T}_{U}^{\beta}(m)$ with
$t_{ii} >0$, $i = 1, \dots,m$) then $(d\mathbf{S})$ (the Lebesgue measure in $\mathfrak{S}_{m}^{\beta}$
or in $\mathfrak{T}_{U}^{\beta}(m)$) denotes the exterior product of the $m(m-1)\beta/2 + m$ functionally
independent variables,
$$
  (d\mathbf{S}) = \bigwedge_{i=1}^{m} ds_{ii}\bigwedge_{i < j}^{m}\bigwedge_{k = 1}^{\beta}
                      ds_{ij}^{(k)}.
$$
Observe, that for the Lebesgue measure $(d\mathbf{S})$ defined thus, it is required that $\mathbf{S} \in
\mathfrak{P}_{m}^{\beta}$, that is, $\mathbf{S}$ must be a non singular Hermitian matrix (Hermitian
definite positive matrix).

If $\boldsymbol{\Lambda} \in \mathfrak{D}_{m}^{\beta}$ then $(d\boldsymbol{\Lambda})$ (the Legesgue
measure in $\mathfrak{D}_{m}^{\beta}$) denotes the exterior product of the $\beta m$ functionally
independent variables

$$
  (d\boldsymbol{\Lambda}) = \bigwedge_{i = 1}^{n}\bigwedge_{k = 1}^{\beta}d\lambda_{i}^{(k)}.
$$
If $\mathbf{H}_{1} \in \mathcal{V}_{m,n}^{\beta}$  is such that $\mathbf{H}_{1} = (\mathbf{h}_{1}, \dots,
\mathbf{h}_{m})$, where $\mathbf{h}_{j}$, $j = 1, \dots,m$ are their columns, then
\begin{equation}\label{eq0}
    (\mathbf{H}^{*}_{1}d\mathbf{H}_{1}) = \bigwedge_{i=1}^{m} \bigwedge_{j =i+1}^{n}
  \mathbf{h}_{j}^{*}d\mathbf{h}_{i},
\end{equation}
where the partitioned matrix $\mathbf{H} = (\mathbf{H}_{1}|\mathbf{H}_{2}) = (\mathbf{h}_{1}, \dots,
\mathbf{h}_{m}|\mathbf{h}_{m+1}, \dots, \mathbf{h}_{n}) \in \mathfrak{U}^{\beta}(n)$, with
$\mathbf{H}_{2} = (\mathbf{h}_{m+1}, \dots, \mathbf{h}_{n})$. It can be proved that this differential
form does not depend on the choice of the $\mathbf{H}_{2}$ matrix. When $n = 1$;
$\mathcal{V}^{\beta}_{m,1}$ defines the unit sphere in $\mathfrak{A}^{m}$. This is, of course, an
$n\beta-1$- dimensional surface in $\mathfrak{A}^{m}$.

The surface area or volume of the Stiefel manifold $\mathcal{V}^{\beta}_{m,n}$ is
\begin{equation}\label{vol}
    \Vol(\mathcal{V}^{\beta}_{m,n}) = \int_{\mathbf{H}_{1} \in
  \mathcal{V}^{\beta}_{m,n}} (\mathbf{H}_{1}d\mathbf{H}^{*}_{1}) =
  \frac{2^{m}\pi^{mn\beta/2}}{\Gamma^{\beta}_{m}[n\beta/2]},
\end{equation}
where $\Gamma^{\beta}_{m}[a]$ denotes the multivariate \emph{Gamma function} for the space
$\mathfrak{S}_{m}^{\beta}$. This can be obtained as a particular case of the \emph{generalised gamma
function of weight $\kappa$} for the space $\mathfrak{S}^{\beta}_{m}$ with $\kappa = (k_{1}, k_{2},
\dots, k_{m}) \in \Re^{m}$, taking $\kappa =(0,0,\dots,0) \in \Re^{m}$ and which for $\re(a) \geq
(m-1)\beta/2 - k_{m}$ is defined by, see \citep{gr:87},
\begin{eqnarray}\label{int1}
  \Gamma_{m}^{\beta}[a,\kappa] &=& \displaystyle\int_{\mathbf{A} \in \mathfrak{P}_{m}^{\beta}}
  \etr\{-\mathbf{A}\} |\mathbf{A}|^{a-(m-1)\beta/2 - 1} q_{\kappa}(\mathbf{A}) (d\mathbf{A}) \\
&=& \pi^{m(m-1)\beta/4}\displaystyle\prod_{i=1}^{m} \Gamma[a + k_{i}
    -(i-1)\beta/2]\nonumber\\ \label{gammagen1}
&=& [a]_{\kappa}^{\beta} \Gamma_{m}^{\beta}[a],
\end{eqnarray}
where $\etr(\cdot) = \exp(\tr(\cdot))$, $|\cdot|$ denotes the determinant, and for $\mathbf{A} \in
\mathfrak{S}_{m}^{\beta}$
\begin{equation}\label{hwv}
    q_{\kappa}(\mathbf{A}) = |\mathbf{A}_{m}|^{k_{m}}\prod_{i = 1}^{m-1}|\mathbf{A}_{i}|^{k_{i}-k_{i+1}}
\end{equation}
with $\mathbf{A}_{p} = (a_{rs})$, $r,s = 1, 2, \dots, p$, $p = 1,2, \dots, m$ is termed the \emph{highest
weight vector}, see \citep{gr:87}. Also,
\begin{eqnarray*}
  \Gamma_{m}^{\beta}[a] &=& \displaystyle\int_{\mathbf{A} \in \mathfrak{P}_{m}^{\beta}}
  \etr\{-\mathbf{A}\} |\mathbf{A}|^{a-(m-1)\beta/2 - 1}(d\mathbf{A}) \\ \label{cgamma}
    &=& \pi^{m(m-1)\beta/4}\displaystyle\prod_{i=1}^{m} \Gamma[a-(i-1)\beta/2],
\end{eqnarray*}
and $\re(a)> (m-1)\beta/2$.

In other branches of mathematics the \textit{highest weight vector} $q_{\kappa}(\mathbf{A})$ is also
termed the \emph{generalised power} of $\mathbf{A}$ and is denoted as $\Delta_{\kappa}(\mathbf{A})$, see
\citep{fk:94} and \citep{hl:01}.

Additional properties of $q_{\kappa}(\mathbf{A})$, which are immediate consequences of the definition of
$q_{\kappa}(\mathbf{A})$ and the following property 1, are:
\begin{enumerate}
  \item Let $\mathbf{A} = \mathbf{L}^{*}\mathbf{DL}$ be the L'DL decomposition of $\mathbf{A} \in \mathfrak{P}_{m}^{\beta}$,
        where $\mathbf{L} \in \mathfrak{T}_{U}^{\beta}(m)$ with $l_{ii} = 1$, $i = 1, 2, \ldots ,m$ and
        $\mathbf{D} = \diag(\lambda_{1}, \dots, \lambda_{m})$, $\lambda_{i} \geq 0$, $i = 1, 2, \ldots
        ,m$. Then
        \begin{equation}\label{qk1}
          q_{\kappa}(\mathbf{A}) = \prod_{i=1}^{m} \lambda_{i}^{k_{i}}.
        \end{equation}
      \item
      \begin{equation}\label{qk2}
        q_{\kappa}(\mathbf{A}^{-1}) =  q_{-\kappa^{*}}^{*}(\mathbf{A}),
      \end{equation}
      where $\kappa^{*}=(k_{m}, k_{m-1}, \dots,k_{1})$, $-\kappa^{*}=(-k_{m}, -k_{m-1},
      \dots,-k_{1})$,
      \begin{equation}\label{hhwv}
         q_{\kappa}^{*}(\mathbf{A}) = |\mathbf{A}_{m}|^{k_{m}}\prod_{i = 1}^{m-1}|\mathbf{A}_{i}|^{k_{i}-k_{i+1}}
      \end{equation}
      and
      \begin{equation}\label{qqk1}
        q_{\kappa}^{*}(\mathbf{A}) = \prod_{i=1}^{m} \lambda_{i}^{k_{m-i+1}},
      \end{equation}
      see \cite[pp. 126-127 and Proposition VII.1.5]{fk:94}.

      Alternatively, let $\mathbf{A} = \mathbf{T}^{*}\mathbf{T}$ the Cholesky's decomposition of
  matrix $\mathbf{A} \in \mathfrak{P}_{m}^{\beta}$, with $\mathbf{T}=(t_{ij}) \in
  \mathfrak{T}_{U}^{\beta}(m)$, then $\lambda_{i} = t_{ii}^{2}$, $t_{ii} \geq 0$, $i = 1, 2,
  \ldots ,m$. See \cite[p. 931, first paragraph]{hl:01}, \cite[p. 390, lines -11 to
  -16]{hlz:05} and \cite[p.5, lines 1-6]{k:14}.
  \item if $\kappa = (p, \dots, p)$, then
    \begin{equation}\label{qk3}
        q_{\kappa}(\mathbf{A}) = |\mathbf{A}|^{p},
    \end{equation}
    in particular if $p=0$, then $q_{\kappa}(\mathbf{A}) = 1$.
  \item if $\tau = (t_{1}, t_{2}, \dots, t_{m})$, $t_{1}\geq t_{2}\geq \cdots \geq t_{m} \geq
  0$, then
    \begin{equation}\label{qk41}
        q_{\kappa+\tau}(\mathbf{A}) = q_{\kappa}(\mathbf{A})q_{\tau}(\mathbf{A}),
    \end{equation}
    in particular if $\tau = (p,p, \dots, p)$,  then
    \begin{equation}\label{qk42}
        q_{\kappa+\tau}(\mathbf{A}) \equiv q_{\kappa+p}(\mathbf{A}) = |\mathbf{A}|^{p} q_{\kappa}(\mathbf{A}).
    \end{equation}
    \item Finally, for $\mathbf{B} \in \mathfrak{T}_{U}^{\beta}(m)$  in such a manner that $\mathbf{C} =
    \mathbf{B}^{*}\mathbf{B} \in \mathfrak{S}_{m}^{\beta}$,
    \begin{equation}\label{qk5}
        q_{\kappa}(\mathbf{B}^{*}\mathbf{AB}) = q_{\kappa}(\mathbf{C})q_{\kappa}(\mathbf{A})
    \end{equation}
    and
    \begin{equation}\label{qk6}
        q_{\kappa}(\mathbf{B}^{*-1}\mathbf{A}\mathbf{B}^{-1}) = (q_{\kappa}(\mathbf{C}))^{-1}q_{\kappa}(\mathbf{A})
        = q_{-\kappa}(\mathbf{C})q_{\kappa}(\mathbf{A}),
    \end{equation}
see \citet[p. 776, eq. (2.1)]{hlz:08}.
\end{enumerate}
\begin{remark}
Let $\mathcal{P}(\mathfrak{S}_{m}^{\beta})$ denote the algebra of all polynomial functions on
$\mathfrak{S}_{m}^{\beta}$, and $\mathcal{P}_{k}(\mathfrak{S}_{m}^{\beta})$ the subspace of homogeneous
polynomials of degree $k$ and let $\mathcal{P}^{\kappa}(\mathfrak{S}_{m}^{\beta})$ be an irreducible
subspace of $\mathcal{P}(\mathfrak{S}_{m}^{\beta})$ such that
$$
  \mathcal{P}_{k}(\mathfrak{S}_{m}^{\beta}) = \sum_{\kappa}\bigoplus
  \mathcal{P}^{\kappa}(\mathfrak{S}_{m}^{\beta}).
$$
Note that $q_{\kappa}$ is a homogeneous polynomial of degree $k$, moreover $q_{\kappa} \in
\mathcal{P}^{\kappa}(\mathfrak{S}_{m}^{\beta})$, see \citep{gr:87}. \qed
\end{remark}
In (\ref{gammagen1}), $[a]_{\kappa}^{\beta}$ denotes the generalised Pochhammer symbol of weight
$\kappa$, defined as
\begin{eqnarray*}
  [a]_{\kappa}^{\beta} &=& \prod_{i = 1}^{m}(a-(i-1)\beta/2)_{k_{i}}\\
    &=& \frac{\pi^{m(m-1)\beta/4} \displaystyle\prod_{i=1}^{m}
    \Gamma[a + k_{i} -(i-1)\beta/2]}{\Gamma_{m}^{\beta}[a]} \\
    &=& \frac{\Gamma_{m}^{\beta}[a,\kappa]}{\Gamma_{m}^{\beta}[a]},
\end{eqnarray*}
where $\re(a) > (m-1)\beta/2 - k_{m}$ and
$$
  (a)_{i} = a (a+1)\cdots(a+i-1),
$$
is the standard Pochhammer symbol.

An alternative definition of the generalised gamma function of weight $\kappa$ is proposed by
\citep{k:66}, which is defined as%
\begin{eqnarray}\label{int2}
  \Gamma_{m}^{\beta}[a,-\kappa] &=& \displaystyle\int_{\mathbf{A} \in \mathfrak{P}_{m}^{\beta}}
    \etr\{-\mathbf{A}\} |\mathbf{A}|^{a-(m-1)\beta/2 - 1} q_{\kappa}(\mathbf{A}^{-1})
    (d\mathbf{A}) \\
&=& \pi^{m(m-1)\beta/4}\displaystyle\prod_{i=1}^{m} \Gamma[a - k_{i}
    -(m-i)\beta/2] \nonumber\\ \label{gammagen2}
&=& \displaystyle\frac{(-1)^{k} \Gamma_{m}^{\beta}[a]}{[-a +(m-1)\beta/2
    +1]_{\kappa}^{\beta}} ,
\end{eqnarray}
where $\re(a) > (m-1)\beta/2 + k_{1}$.

In addition, observe that by (\ref{eq0}) and (\ref{vol})
$$
  (d\mathbf{H}_{1}) = \frac{1}{\Vol\left(\mathcal{V}^{\beta}_{m,n}\right)}
    (\mathbf{H}_{1}^{*}d\mathbf{H}_{1}) = \frac{\Gamma^{\beta}_{m}[n\beta/2]}{2^{m}
    \pi^{mn\beta/2}}(\mathbf{H}_{1}^{*}d\mathbf{H}_{1}).
$$
is the \emph{normalised invariant measure on} $\mathcal{V}^{\beta}_{m,n}$ and $(d\mathbf{H})$, i.e., with
$(m = n)$, it defines the \emph{normalised Haar measure} on $\mathfrak{U}^{\beta}(m)$.

Finally, the ${}_{p}F_{q}^{\beta}$ is the generalised hypergeometric function with matrix argument,
defined as
$$
  {}_{p}F_{q}^{\beta}(a_{1}, \dots,a_{p};b_{1}, \dots,b_{q}; \mathbf{X}) = \sum_{t =0}^{\infty}\sum_{\tau}
  \frac{[a_{1}]_{\tau}^{\beta} \cdots [a_{p}]_{\tau}^{\beta}}{[b_{1}]_{\tau}^{\beta} \cdots [b_{q}]_{\tau}^{\beta}}
  \frac{C_{\tau}^{\beta}(\mathbf{X})}{t!}
$$
where, $\sum_{\tau}$, denotes summation over all partition $\tau = (t_{1}, t_{2},\dots, t_{m})$,
$t_{1}\geq t_{2}\geq \cdots t_{m} \geq 0$, of $t$, $\sum_{i=1}^{m} t_{i}=t$ and $t_{1}, t_{2},\dots,
t_{m}$ are nonnegative integers, and $C_{\tau}^{\beta}(\mathbf{X})$ is the Jack polynomial of $\mathbf{X}
\in {\mathfrak S}_{m}^{\beta}$ corresponding to $t$, see \citep{gr:87} and \citep{dg:14}.

Jack polynomials for real normed division algebras are also termed spherical functions of symmetric cones
in the abstract algebra context, see \citep{S:97}. In addition, in the statistical literature, they are
termed real, complex, quaternion and octonion zonal polynomials, or, generically, \textit{general zonal
polynomials}, see \citep{j:64}, Chapter 7 in \citep{m:82}, \citep{k:84} and \citep{lx:09}. This section
is completed, remembering the following result: from \citep[Equation 4.8(2) and Definition 5.3]{gr:87} is
obtained
\begin{equation}\label{jpq}
    C_{\kappa}^{\beta}(\mathbf{X}) = C_{\kappa}^{\beta}(\mathbf{I}_{m})
    \int_{\mathbf{H} \in \mathfrak{U}^{\beta}(m)} q_{\kappa}(\mathbf{H}^{*}\mathbf{XH})(d\mathbf{H})
\end{equation}
for all $\mathbf{X} \in \mathfrak{S}_{m}^{\beta}$; where $(d\mathbf{H})$ is the normalised Haar measure
on $\mathfrak{U}^{\beta}(m)$, also see \citep{dg:14}.

\section{Integration}\label{sec3}
First consider the following result proposed by \citep{zf:90} for real case and extended to real normed
division algebras by \citep{dg:14}.
\begin{lemma}\label{cfu}
The characteristic function of $\mathbf{H}_{1} \in \mathcal{V}^{\beta}_{m,n}$, the normalised invariant
measure on $\mathcal{V}^{\beta}_{m,n}$ is
\begin{eqnarray}
  \phi_{_{\mathbf{H}_{1}}}(\mathbf{T}) &=& \int_{\mathbf{H}_{1} \in \mathcal{V}^{\beta}_{m,n}}
   \etr\{i\mathbf{H}_{1}\mathbf{T}^{*}\}(d\mathbf{H}_{1}) \nonumber\\
   &=& {}_{0}F_{1}^{\beta}(\beta n/2, -\mathbf{TT}^{*}/4) \nonumber \\ \label{chu}
   &=&\displaystyle \sum_{t =0}^{\infty}\sum_{\tau}\frac{1}{[\beta
   n/2]_{\tau}^{\beta}}\frac{C_{\tau}^{\beta}(-\mathbf{TT}^{*}/4)}{t!},
\end{eqnarray}
where, $\sum_{\tau}$, denotes summation over all partition $\tau = (t_{1}, t_{2},\dots, t_{m})$,
$t_{1}\geq t_{2}\geq \cdots t_{m} \geq 0$, of $t$, $\sum_{i=1}^{m} t_{i}=t$ and $t_{1}, t_{2},\dots,
t_{m}$ are nonnegative integers.
\end{lemma}

The original version of the following result was obtained for real case in \citep{xf:90}. Next, we set
the version of this result for real normed division algebras.

First consider the following concept: let's use the complexification $\mathfrak{S}_{m}^{\beta,
\mathfrak{C}} = \mathfrak{S}_{m}^{\beta} + i \mathfrak{S}_{m}^{\beta}$ of $\mathfrak{S}_{m}^{\beta}$.
That is, $\mathfrak{S}_{m}^{\beta, \mathfrak{C}}$ consist of all matrices $\mathbf{Z} \in
(\mathfrak{F^{\mathfrak{C}}})^{m \times m}$ of the form $\mathbf{Z} = \mathbf{X} + i\mathbf{Y}$, with
$\mathbf{X}, \mathbf{Y} \in \mathfrak{S}_{m}^{\beta}$. It comes to $\mathbf{X} = \re(\mathbf{Z})$ and
$\mathbf{Y} = \im(\mathbf{Z})$ as the \emph{real and imaginary parts} of $\mathbf{Z}$, respectively. The
\emph{generalised right half-plane} $\boldsymbol{\Phi}_{m}^{\beta} = \mathfrak{P}_{m}^{\beta} + i
\mathfrak{S}_{m}^{\beta}$ in $\mathfrak{S}_{m}^{\beta,\mathfrak{C}}$ consists of all $\mathbf{Z} \in
\mathfrak{S}_{m}^{\beta,\mathfrak{C}}$ such that $\re(\mathbf{Z}) \in \mathfrak{P}_{m}^{\beta}$, see
\citep[p. 801]{gr:87}.

\begin{lemma}\label{teo2i}
Let $\mathbf{Z} \in \boldsymbol{\Phi}_{m}^{\beta}$ and $\mathbf{U} \in \mathfrak{S}_{m}^{\beta}$. Then
\begin{eqnarray}\label{runze22}
    \int_{\mathbf{X} \in \mathfrak{P}_{m}^{\beta}} f(\mathbf{Z}^{1/2}\mathbf{XZ}^{1/2}) |\mathbf{X}|^{a -(m-1)
    \beta/2-1} C_{\tau}^{\beta}\left(\mathbf{X}\mathbf{U} \right)
    (d\mathbf{X}) \hspace{3cm} \nonumber\\\label{jpe2}
    =\displaystyle\frac{J(\mathbf{I}_{m})}{C_{\tau}^{\beta}(\mathbf{I}_{m})}
     |\mathbf{Z}|^{-a} C_{\tau}^{\beta}(\mathbf{UZ}^{-1}),
\end{eqnarray}
where $\mathbf{Z}^{1/2}$ is the positive definite square root of $\mathbf{Z}$, i.e.
$\mathbf{Z}^{1/2}\mathbf{Z}^{1/2} = \mathbf{Z}$, $\re(a)> (m-1)\beta/2 - t_{m}$, $\tau = (t_{1},
t_{2},\dots, t_{m})$, $t_{1}\geq t_{2}\geq \cdots t_{m} \geq 0$, $\sum_{i=1}^{m} t_{i}=t$ and $t_{1},
t_{2},\dots, t_{m}$ are nonnegative integers, and
$$
  J(\mathbf{I}_{m}) = \int_{\mathbf{X} \in \mathfrak{P}_{m}^{\beta}} f(\mathbf{X}) |\mathbf{X}|^{a -(m-1)
    \beta/2-1} C_{\tau}^{\beta}\left(\mathbf{X}\right)
    (d\mathbf{X}).
$$
\end{lemma}
\begin{proof}
The proof is a verbatim copy of given by \citep{xf:90}, only take notice that
\begin{enumerate}
  \item If $\mathbf{Y} = \mathbf{H}^{*}\mathbf{X}\mathbf{H}$ for $\mathbf{H} \in
  \mathfrak{U}^{\beta}(m)$, then $(d\mathbf{Y}) = (d\mathbf{X})$, and
  \item if $\mathbf{B} = \mathbf{Z}^{1/2}\mathbf{XZ}^{1/2}$ with $\mathbf{B}, \mathbf{X} \in
  \mathfrak{P}_{m}^{\beta}$, then by \citep[Proposition 2]{dggj:13}, $(d\mathbf{B}) =
  |\mathbf{Z}|^{(m-1)\beta/2+1}(d\mathbf{X})$, where $\mathbf{Z}^{1/2}$ is the positive definite square root of $\mathbf{Z}$,
  such that $\mathbf{Z}^{1/2}\mathbf{Z}^{1/2} = \mathbf{Z}$.
\end{enumerate}
\end{proof}

The following result shows the extension of the Wishart's integral for real normed division algebras
stated in \citep{dg:13}.

\begin{lemma}\label{teo3i}
Let $\mathbf{Y} \in \mathfrak{L}^{\beta}_{n,m}$.
$$
   \int_{\mathbf{Y}^{*}\mathbf{Y}= \mathbf{R}} f(\mathbf{Y}^{*}\mathbf{Y}) d(\mathbf{Y}) =
   \frac{\pi^{\beta mn/2}}{\Gamma_{m}^{\beta}[\beta n/2]}|\mathbf{R}|^{\beta(n - m + 1)/2 - 1}
   f(\mathbf{R}).
$$
where $\re(\beta n/2)> (m-1)\beta/2$.
\end{lemma}

Finally, the Theorem 1 in \citep{l:93} is generalised for real normed division algebras. But prior
consider the following definition, see \citep{dggj:13}.

\begin{definition}\label{defLS}
Let $\mathbf{X} \in \mathfrak{L}_{n,m}^{\beta}$ be a random matrix, then if $\mathbf{X} \build{=}{d}{}
\boldsymbol{\Xi}\mathbf{X}$ for every $\boldsymbol{\Xi} \in \mathfrak{U}^{\beta}(n)$, $\mathbf{X}$ is
termed left-spherical. Where $\build{=}{d}{}$ signifies that the two sides have the same distribution.
\end{definition}

In addition, note that if $\mathbf{X}$ has density, respect to the Lebesgue measure, it has the form
$f(\beta \mathbf{X}^{*}\mathbf{X})$ and $\mathbf{X}$ can be factorised as
\begin{equation}\label{sr}
    \mathbf{X} = \mathbf{H}_{1}\mathbf{A}
\end{equation}
where $\mathbf{A}_{m \times m}$  is not unique, and it is independent of $\mathbf{H}_{1}$ and
$\mathbf{H}_{1}$ has a normalised invariant measure on $\mathcal{V}^{\beta}_{m,n}$.
\begin{theorem}\label{teocf}
Suppose that $\mathbf{X}$ is as in Definition \ref{defLS}. Then the characteristic function of
$\mathbf{X}$ can be expressed as
$$
  \phi_{_{\mathbf{X}}}(\mathbf{T})=  \sum_{t =0}^{\infty}\sum_{\tau}\frac{C_{\tau}^{\beta}(-\mathbf{TT}^{*}/4)}
  {[\beta n/2]_{\tau}^{\beta}\ C_{\tau}^{\beta}(\mathbf{I}_{m})\ t!} \E\left(C_{\tau}^{\beta}(\mathbf{R})\right)
$$
where $\mathbf{R} = \mathbf{X}^{*}\mathbf{X}$, $\re(\beta n/2)> (m-1)\beta/2 - t_{m}$, $\tau = (t_{1},
t_{2},\dots, t_{m})$, $t_{1}\geq t_{2}\geq \cdots t_{m} \geq 0$, $\sum_{i=1}^{m} t_{i}=t$ and $t_{1},
t_{2},\dots, t_{m}$ are nonnegative integers.
\end{theorem}
\begin{proof}By (\ref{sr}) and Lemma \ref{cfu},
\begin{eqnarray*}
   \phi_{_{\mathbf{X}}}(\mathbf{T})&=& \E_{_{\mathbf{X}}}(\etr\{i\mathbf{XT}^{*}\})
   = \E_{_{\mathbf{H}_{1},\mathbf{A}}}(\etr\{i\mathbf{H}_{1}\mathbf{A}\mathbf{T}^{*}\})\\
   &=& \E_{_{\mathbf{H}_{1},\mathbf{A}}}(\etr\{i\mathbf{H}_{1}(\mathbf{T}\mathbf{A}^{*})^{*}\}) =
   \E_{_{\mathbf{A}}} \left( \phi_{_{\mathbf{H}_{1}}}(\mathbf{T}\mathbf{A}^{*})\right)\\
   && \mbox{because $\mathbf{H}_{1}$ and $\mathbf{A}$ are independent,}\\
   &=& \E_{_{\mathbf{A}}} \left({}_{0}F_{1}^{\beta}(\beta n/2,
   -\mathbf{T}\mathbf{A}^{*}(\mathbf{T}\mathbf{A}^{*})^{*}/4)\right)\\
   &=& \E_{_{\mathbf{A}}} \left({}_{0}F_{1}^{\beta}(\beta n/2,
   -\mathbf{T}^{*}\mathbf{T}\mathbf{A}^{*}\mathbf{A}/4)\right)\\
   &=& \E_{_{\mathbf{R}}} \left({}_{0}F_{1}^{\beta}(\beta n/2,
   -\mathbf{T}^{*}\mathbf{T}\mathbf{R}/4)\right) \quad \mbox{with $\mathbf{R} =
   \mathbf{A}^{*}\mathbf{A}$}\\
   && \mbox{thus, by Lemma \ref{cfu} and Lemma \ref{teo3i},}\\
   &=& \displaystyle \sum_{t =0}^{\infty}\sum_{\tau}\frac{\pi^{\beta mn/2}}{[\beta
   n/2]_{\tau}^{\beta}\Gamma_{m}^{\beta}[\beta n/2]t!}\\
   && \hspace{1.5cm} \times \int_{\mathbf{R} \in \mathfrak{P}_{m}^{\beta}}|\mathbf{R}|^{\beta(n - m + 1)/2 - 1}f(\mathbf{R})
   C_{\tau}^{\beta}(-\mathbf{T}^{*}\mathbf{T}\mathbf{R}/4)(d\mathbf{R})\\
   && \mbox{hence, by Lemma \ref{teo2i},}\\
   &=& \displaystyle \sum_{t =0}^{\infty}\sum_{\tau}\frac{C_{\tau}^{\beta}(-\mathbf{T}^{*}\mathbf{T}/4)}{[\beta
   n/2]_{\tau}^{\beta}C_{\tau}^{\beta}(\mathbf{I}_{m})t!}\\
   && \hspace{1.5cm} \times \int_{\mathbf{R} \in \mathfrak{P}_{m}^{\beta}}\frac{\pi^{\beta mn/2}}{\Gamma_{m}^{\beta}[\beta n/2]}
   |\mathbf{R}|^{\beta(n - m + 1)/2 - 1}f(\mathbf{R}) C_{\tau}^{\beta}(\mathbf{R})(d\mathbf{R})\\
   &=& \displaystyle \sum_{t =0}^{\infty}\sum_{\tau}\frac{C_{\tau}^{\beta}(-\mathbf{T}^{*}\mathbf{T}/4)}{[\beta
   n/2]_{\tau}^{\beta}\ C_{\tau}^{\beta}(\mathbf{I}_{m})\ t!}
   \E_{_{\mathbf{R}}}\left( C_{\tau}^{\beta}(\mathbf{R})\right).
\end{eqnarray*}
\end{proof}

\section{Kotz-Riesz distributions}\label{sec4}

This section introduce two versions of the Kotz-Riesz distributions. With this purpose in mind, first are
defined the spherical Kotz-Riesz distributions.

\begin{definition}\label{defnSKR}
Let $\kappa = (k_{1}, k_{2}, \dots, k_{m}) \in \Re^{m}$.
\begin{enumerate}
  \item Then it is said that $\mathbf{Z}$ has a spherical Kotz-Riesz distribution of type I if its density function is
  \begin{equation}\label{dfSKR1}
    \frac{\beta^{mn\beta/2}\Gamma_{m}^{\beta}[n\beta/2]}{\pi^{mn\beta/2}\Gamma_{m}^{\beta}[n\beta/2,\kappa]}
    \etr\left\{- \beta\tr \mathbf{Z}^{*} \mathbf{Z}\right\}
   q_{\kappa}\left [\beta \mathbf{Z}^{*}\mathbf{Z}\right ](d\mathbf{Z})
  \end{equation}
  for $\mathbf{Z} \in \mathfrak{L}^{\beta}_{n,m}$, $\re(n\beta/2) > (m-1)\beta/2 - k_{m}$;  denoting this fact as
  $$
    \mathbf{Z} \sim \mathcal{K}\mathfrak{R}^{\beta, I}_{n \times m}
    (\kappa,\boldsymbol{0}, \mathbf{I}_{n}, \mathbf{I}_{m}).
  $$
  \item Then it is said that $\mathbf{Z}$ has a Kotz-Riesz distribution of type II if its density function is
  \begin{equation}\label{dfSKR2}
  \frac{\beta^{mn\beta/2}\Gamma_{m}^{\beta}[n\beta/2]}{\pi^{mn\beta/2}\Gamma_{m}^{\beta}[n\beta/2,\kappa]}
    \etr\left\{- \beta\tr \mathbf{Z}^{*} \mathbf{Z}\right\}
   q_{\kappa}\left [(\beta \mathbf{Z}^{*}\mathbf{Z})^{-1}\right ](d\mathbf{Z})
  \end{equation}
  for $\mathbf{Z} \in \mathfrak{L}^{\beta}_{n,m}$, $\re(n\beta/2) > (m-1)\beta/2 + k_{1}$;  denoting this fact as
  $$
    \mathbf{Z} \sim \mathcal{K}\mathfrak{R}^{\beta, II}_{n \times m}
    (\kappa,\boldsymbol{0}, \mathbf{I}_{n}, \mathbf{I}_{m}).
  $$
\end{enumerate}
\end{definition}
Now the \emph{elliptical Kotz-Riesz distribution} or simply termed \emph{Kotz-Riesz distribution} is
obtained.
\begin{theorem}\label{teoEKR}
Let $\boldsymbol{\Sigma} \in \boldsymbol{\Phi}_{m}^{\beta}$, $\boldsymbol{\Theta} \in
\boldsymbol{\Phi}_{n}^{\beta}$, $\boldsymbol{\mu} \in \mathfrak{L}^{\beta}_{n,m}$ and  $\kappa = (k_{1},
k_{2}, \dots, k_{m}) \in \Re^{m}$. And let $\mathbf{Z} \sim \mathcal{K}\mathfrak{R}^{\beta, I}_{n \times
m}(\kappa,\boldsymbol{0}, \mathbf{I}_{n}, \mathbf{I}_{m})$. Also, defines $\mathbf{X} = \boldsymbol{\mu}+
\u(\boldsymbol{\Theta})^{*}\mathbf{Z}\u(\boldsymbol{\Sigma})$, where $\u(\mathbf{B}) \in
\mathfrak{T}_{U}^{\beta}(m)$ is such that $\mathbf{B} = \u(\mathbf{B})^{*}\u(\mathbf{B})$ is the Cholesky
decomposition of $\mathbf{B}$.
\begin{enumerate}
  \item Then it is said that $\mathbf{X}$ has a Kotz-Riesz distribution of type I and its density function is
  $$
    \frac{\beta^{mn\beta/2+\sum_{i = 1}^{m}k_{i}}\Gamma_{m}^{\beta}[n\beta/2]}{\pi^{mn\beta/2}\Gamma_{m}^{\beta}[n\beta/2,\kappa]
    |\boldsymbol{\Sigma}|^{n\beta/2}|\boldsymbol{\Theta}|^{m\beta/2}}\hspace{4cm}
  $$
  $$\hspace{1cm}
    \times \etr\left\{- \beta\tr \left [\boldsymbol{\Sigma}^{-1} (\mathbf{X} - \boldsymbol{\mu})^{*}
    \boldsymbol{\Theta}^{-1}(\mathbf{X} - \boldsymbol{\mu})\right ]\right\}
  $$
  \begin{equation}\label{dfEKR1}\hspace{3.1cm}
    \times q_{\kappa}\left [\u(\boldsymbol{\Sigma})^{*-1} (\mathbf{X} - \boldsymbol{\mu})^{*}
    \boldsymbol{\Theta}^{-1}(\mathbf{X} - \boldsymbol{\mu})\u(\boldsymbol{\Sigma})^{-1}\right ](d\mathbf{X})
  \end{equation}
  for $\mathbf{X} \in \mathfrak{L}^{\beta}_{n,m}$, $\re(n\beta/2) > (m-1)\beta/2 - k_{m}$;  denoting this fact as
  $$
    \mathbf{X} \sim \mathcal{K}\mathfrak{R}^{\beta, I}_{n \times m}
    (\kappa,\boldsymbol{\mu}, \boldsymbol{\Theta}, \boldsymbol{\Sigma}).
  $$
  \item Then it is said that $\mathbf{X}$ has a Kotz-Riesz distribution of type II and its density function is
  $$
    \frac{\beta^{mn\beta/2-\sum_{i = 1}^{m}k_{i}}\Gamma_{m}^{\beta}[n\beta/2]}{\pi^{mn\beta/2}\Gamma_{m}^{\beta}[n\beta/2,-\kappa]
    |\boldsymbol{\Sigma}|^{n\beta/2}|\boldsymbol{\Theta}|^{m\beta/2}}\hspace{4cm}
  $$
  $$
    \times \etr\left\{- \beta\tr \left [\boldsymbol{\Sigma}^{-1} (\mathbf{X} - \boldsymbol{\mu})^{*}
    \boldsymbol{\Theta}^{-1}(\mathbf{X} - \boldsymbol{\mu})\right ]\right\}
  $$
  \begin{equation}\label{dfEKR2}\hspace{2.5cm}
    \times q_{\kappa}\left [\left(\u(\boldsymbol{\Sigma})^{*-1} (\mathbf{X} - \boldsymbol{\mu})^{*}
    \boldsymbol{\Theta}^{-1}(\mathbf{X} - \boldsymbol{\mu})\u(\boldsymbol{\Sigma})^{-1/2}\right)^{-1}\right ](d\mathbf{X})
  \end{equation}
  for $\mathbf{X} \in \mathfrak{L}^{\beta}_{n,m}$, $\re(n\beta/2) > (m-1)\beta/2 + k_{1}$;  denoting this fact as
  $$
    \mathbf{X} \sim \mathcal{K}\mathfrak{R}^{\beta, II}_{n \times m}
    (\kappa,\boldsymbol{\mu}, \boldsymbol{\Theta}, \boldsymbol{\Sigma}).
  $$
\end{enumerate}
\end{theorem}
\begin{proof}
This is an immediate consequence of Theorem 3 in \citep{dggj:13} and the fact is that, for $a$ scalar
$$
  q_{\kappa}(a\mathbf{A}) = q_{\kappa}(a\mathbf{I}_{m}) q_{\kappa}(\mathbf{A})=
  a^{\sum_{i=1}^{m}k_{i}}q_{\kappa}(\mathbf{A}).
$$
\end{proof}

Observe that, if $\kappa = (0,0, \dots,0)$  and $\boldsymbol{\Sigma} = 2\boldsymbol{\Sigma}$ in two
densities in Theorem \ref{teoEKR} the matrix multivariate normal distribution for real normed division
algebras is obtained, see \citep{dggj:11}. Also, when $\kappa = (l,l, \dots,l)$, $l = a-(m-1)\beta/2 - 1$
and $\boldsymbol{\Sigma} = \boldsymbol{\Sigma}/r$, $r > 1$, the original Kotz type distribution are
obtained i.e. for $s = 1$, in notation of \citep{l:93}.

Next, the main consequence of this article is is achieved by finding the two versions of the Riesz
distributions in terms of Kotz-Riesz distributions.

\begin{theorem}\label{teoRd}
\begin{enumerate}
    \item Assume that $\mathbf{X} \sim \mathcal{K}\mathfrak{R}^{\beta, I}_{n \times m}(\kappa,\boldsymbol{0},
    \boldsymbol{\Theta}, \boldsymbol{\Sigma})$, and define $\mathbf{Y} = \mathbf{X}^{*}\boldsymbol{\Theta}^{-1}\mathbf{X}$. Then
    its density function is
  \begin{equation}\label{dfR1}
    \frac{\beta^{am+\sum_{i = 1}^{m}k_{i}}}{\Gamma_{m}^{\beta}[a,\kappa] |\boldsymbol{\Sigma}|^{a}q_{\kappa}(\boldsymbol{\Sigma})}
    \etr\{-\beta \boldsymbol{\Sigma}^{-1}\mathbf{Y}\}|\mathbf{Y}|^{a-(m-1)\beta/2 - 1}
    q_{\kappa}(\mathbf{Y})(d\mathbf{Y})
  \end{equation}
  for $\mathbf{Y} \in \mathfrak{P}_{m}^{\beta}$ and $\re(a) > (m-1)\beta/2 - k_{m}$; denoting this fact as
  $\mathbf{Y} \sim \mathfrak{R}^{\beta, I}_{m}(a,\kappa, \boldsymbol{\Sigma})$.
  \item Suppose that $\mathbf{X} \sim \mathcal{K}\mathfrak{R}^{\beta, II}_{n \times m}(\kappa,\boldsymbol{0},
    \boldsymbol{\Theta}, \boldsymbol{\Sigma})$, and define $\mathbf{Y} = \mathbf{X}^{*}\mathbf{X}$. Then
    its density function is
  \begin{equation}\label{dfR2}
     \frac{\beta^{am-\sum_{i = 1}^{m}k_{i}}}{\Gamma_{m}^{\beta}[a,-\kappa]
   |\boldsymbol{\Sigma}|^{a}q_{\kappa}(\boldsymbol{\Sigma}^{-1}) \ }\etr\{-\beta \boldsymbol{\Sigma}^{-1}\mathbf{Y}\}
  |\mathbf{Y}|^{a-(m-1)\beta/2 - 1} q_{\kappa}(\mathbf{Y}^{-1}) (d\mathbf{Y})
  \end{equation}
  for $\mathbf{Y} \in \mathfrak{P}_{m}^{\beta}$ and $\re(a) > (m-1)\beta/2 + k_{1}$; denoting this fact as
  $\mathbf{Y} \sim \mathfrak{R}^{\beta, II}_{m}(a,\kappa, \boldsymbol{\Sigma})$.
\end{enumerate}
Where $\kappa = (k_{1}, k_{2}, \dots, k_{m}) \in \Re^{m}$.
\end{theorem}
\begin{proof}
\begin{enumerate}
  \item By applying (\ref{qk5}), the density of $\mathbf{X}$ is
  $$
    \frac{\beta^{mn\beta/2+\sum_{i = 1}^{m}k_{i}}\Gamma_{m}^{\beta}[n\beta/2]}{\pi^{mn\beta/2}\Gamma_{m}^{\beta}[n\beta/2,\kappa]
    |\boldsymbol{\Sigma}|^{n\beta/2}q_{\kappa}(\boldsymbol{\Sigma})|\boldsymbol{\Theta}|^{m\beta/2}}\hspace{4cm}
  $$
  $$\hspace{1cm}
    \times \etr\left\{- \beta\tr \boldsymbol{\Sigma}^{-1} \mathbf{X}^{*}
    \boldsymbol{\Theta}^{-1}\mathbf{X}\right\} q_{\kappa}\left [\mathbf{X}^{*}
    \boldsymbol{\Theta}^{-1}\mathbf{X}\right ](d\mathbf{X}).
  $$
  Now, let $\boldsymbol{\Theta}^{1/2}$ the positive definite square root of $\boldsymbol{\Theta}$ and define
  $\mathbf{V} = \boldsymbol{\Theta}^{-1/2} \mathbf{X}$. Hence, $\mathbf{V} \sim \mathcal{K}\mathfrak{R}^{\beta, I}_{n \times m}
  (\kappa,\boldsymbol{0}, \mathbf{I}_{n}, \boldsymbol{\Sigma})$; moreover, its density is
  $$
    \frac{\beta^{mn\beta/2+\sum_{i = 1}^{m}k_{i}}\Gamma_{m}^{\beta}[n\beta/2]}{\pi^{mn\beta/2}\Gamma_{m}^{\beta}[n\beta/2,\kappa]
    |\boldsymbol{\Sigma}|^{n\beta/2}q_{\kappa}(\boldsymbol{\Sigma})}
     \times \etr\left\{- \beta\tr \boldsymbol{\Sigma}^{-1} \mathbf{V}^{*}
    \mathbf{V}\right\} q_{\kappa}\left [\mathbf{V}^{*}
    \mathbf{V}\right ](d\mathbf{V}).
  $$
  Finally, defining $\mathbf{Y} = \mathbf{V}^{*}\mathbf{V}$ the desired result is obtained by
  applying the Lemma \ref{teo3i} and denoting $a = n\beta/2$.
  \item  It is obtained in analogous way to 1.
\end{enumerate}
\end{proof}

Note that (\ref{dfR1}) and (\ref{dfR2}) are the density functions of Riesz distributions type I and II,
respectively; which were obtained on a special case of the  \emph{Riesz measure} by \citep{fk:94},
\citep{hl:01} and \citep{dg:15a}.

Below, are found the characteristic functions of the Kotz-Riesz distribution type I.  But before, observe
that in order to apply  the Theorem \ref{teocf} it is necessary  to find the expected value of
$C_{\tau}^{\beta}(\mathbf{AY})$ where $\mathbf{Y} = \mathbf{X}^{*}\mathbf{X}$ has a Riesz distribution.
\begin{theorem}\label{EJP}
Assume that $\mathbf{Z} \sim \mathcal{K}\mathfrak{R}^{\beta, I}_{n \times m} (\kappa,\boldsymbol{0},
\mathbf{I}_{n}, \mathbf{I}_{m})$, then
$$
  \E(C_{\tau}^{\beta}(\mathbf{AY}))= \frac{\beta^{mn\beta/2+\sum_{i = 1}^{m}k_{i}}\Gamma_{m}^{\beta}[n\beta/2,\tau+\kappa]}
  {\Gamma_{m}^{\beta}[n\beta/2,\kappa]} C_{\tau}^{\beta}(\mathbf{A})
$$
where $\mathbf{Y} \build{=}{d}{}\mathbf{Z}^{*}\mathbf{Z}$; $\tau = (t_{1}, \dots, t_{m})$, $t_{1}\geq
t_{2}\geq \cdots \geq t_{m} \geq 0$, $t_{1}, t_{2},\dots, t_{m}$ are nonnegative integers; and $\kappa =
(k_{1}, k_{2}, \dots, k_{m})$, $k_{1}\geq k_{2}\geq \cdots \geq k_{m} \geq 0$, $k_{1}, k_{2},\dots,
k_{m}$ are nonnegative integers.
\end{theorem}
\begin{proof}
By (\ref{dfR1}) it is had that $\mathbf{Y} = \mathbf{Z}^{*}\mathbf{Z} \sim \mathfrak{R}^{\beta,
I}_{m}(\beta n/2,\kappa, \mathbf{I}_{m})$. Then $\E(C_{\tau}^{\beta}(\mathbf{AY}))$ is
\begin{eqnarray*}
  &=& \frac{\beta^{am+\sum_{i = 1}^{m}k_{i}}}{\Gamma_{m}^{\beta}[\beta n/2,\kappa]}
    \int_{\mathbf{Y} \in \mathfrak{P}_{m}^{\beta}}C_{\tau}^{\beta}(\mathbf{AY}) \etr\{-\beta \mathbf{Y}\}
    |\mathbf{Y}|^{(n-m-1)\beta/2 - 1} q_{\kappa}(\mathbf{Y})(d\mathbf{Y}) \\
  && \mbox{by Lemma \ref{teo2i}} \\
  &=& \frac{\beta^{am+\sum_{i = 1}^{m}k_{i}}}{\Gamma_{m}^{\beta}[\beta n/2,\kappa]}
  \frac{J(\mathbf{I}_{m})}{C_{\tau}^{\beta}(\mathbf{I}_{m})}C_{\tau}^{\beta}(\mathbf{A}),
\end{eqnarray*}
where by (\ref{jpq})
\begin{eqnarray*}
  J(\mathbf{I}_{m}) &=& \int_{\mathbf{Y} \in \mathfrak{P}_{m}^{\beta}} \etr\{-\beta \mathbf{Y}\}
    |\mathbf{Y}|^{(n-m-1)\beta/2 - 1} q_{\kappa}(\mathbf{Y})C_{\tau}^{\beta}(\mathbf{Y})(d\mathbf{Y}) \\
   &=&\int_{\mathbf{Y} \in \mathfrak{P}_{m}^{\beta}} \etr\{-\beta \mathbf{Y}\}
    |\mathbf{Y}|^{(n-m-1)\beta/2 - 1} q_{\kappa}(\mathbf{Y}) \\
   &&  \hspace{3cm}\times C_{\tau}^{\beta}(\mathbf{I}_{m}) \int_{\mathbf{H} \in \mathfrak{U}^{\beta}(m)}
    q_{\tau}(\mathbf{H}^{*}\mathbf{Y}\mathbf{H})(d\mathbf{H})(d\mathbf{Y}).
\end{eqnarray*}
Thus, proceeding as in the proof of Theorem 5.9 in \cite[p. 802]{gr:87}, and by applying (\ref{qk41})
\begin{eqnarray*}
  J(\mathbf{I}_{m}) &=& C_{\tau}^{\beta}(\mathbf{I}_{m})\int_{\mathbf{Y} \in \mathfrak{P}_{m}^{\beta}} \etr\{-\beta \mathbf{Y}\}
    |\mathbf{Y}|^{(n-m-1)\beta/2 - 1} q_{\kappa}(\mathbf{Y}) q_{\tau}(\mathbf{Y})(d\mathbf{Y})\\
  &=& C_{\tau}^{\beta}(\mathbf{I}_{m})\int_{\mathbf{Y} \in \mathfrak{P}_{m}^{\beta}} \etr\{-\beta \mathbf{Y}\}
    |\mathbf{Y}|^{(n-m-1)\beta/2 - 1} q_{\kappa+\tau}(\mathbf{Y}) (d\mathbf{Y})\\
  &=& C_{\tau}^{\beta}(\mathbf{I}_{m}) \Gamma_{m}^{\beta}[\beta n/2,\kappa+\tau].
\end{eqnarray*}
From where the desired result is achieved.
\end{proof}

The following result gives the characteristic function of the spherical Kotz-Riesz distribution type I.

\begin{theorem}\label{teocfKR1}
Assume that $\mathbf{Z} \sim \mathcal{K}\mathfrak{R}^{\beta, I}_{n \times m} (\kappa,\boldsymbol{0},
\mathbf{I}_{n}, \mathbf{I}_{m})$. Then the characteristic function of $\mathbf{Z}$ can be expressed as
$$
  \phi_{_{\mathbf{Z}}}(\mathbf{T})=  \frac{\beta^{mn\beta/2+\sum_{i = 1}^{m}k_{i}}}
  {\Gamma_{m}^{\beta}[n\beta/2,\kappa]}\sum_{t =0}^{\infty}\sum_{\tau}\frac{\Gamma_{m}^{\beta}[n\beta/2,\tau+\kappa]}
  {[\beta n/2]_{\tau}^{\beta}}\frac{C_{\tau}^{\beta}(-\mathbf{TT}^{*}/4)}{t!},
$$
where, $\re(\beta n/2)> (m-1)\beta/2 - t_{m}$, $\tau = (t_{1}, t_{2},\dots, t_{m})$, $t_{1}\geq t_{2}\geq
\cdots t_{m} \geq 0$, $\sum_{i=1}^{m} t_{i}=t$ and $t_{1}, t_{2},\dots, t_{m}$ are nonnegative integers;
and $\kappa = (k_{1}, k_{2}, \dots, k_{m})$, $k_{1}\geq k_{2}\geq \cdots \geq k_{m} \geq 0$, $k_{1},
k_{2},\dots, k_{m}$ are nonnegative integers.
\end{theorem}
\begin{proof}
This is immediately from Theorem \ref{teocf} and Theorem \ref{EJP}.
\end{proof}

Finally, the characteristic function of the elliptical Kotz-Riesz distribution type I is given.

\begin{corollary}
Assume that $\mathbf{X} \sim \mathcal{K}\mathfrak{R}^{\beta, I}_{n \times m} (\kappa,\boldsymbol{\mu},
\boldsymbol{\Theta}, \boldsymbol{\Sigma})$. Then the characteristic function of $\mathbf{X}$ can be
expressed as
$$
  \phi_{_{\mathbf{X}}}(\mathbf{T})=  \frac{\beta^{mn\beta/2+\sum_{i = 1}^{m}k_{i}}\etr\{i\boldsymbol{\mu}\mathbf{T}^{*}\}}
  {\Gamma_{m}^{\beta}[n\beta/2,\kappa]}\sum_{t =0}^{\infty}\sum_{\tau}\frac{\Gamma_{m}^{\beta}[n\beta/2,\tau+\kappa]}
  {[\beta n/2]_{\tau}^{\beta}}\frac{C_{\tau}^{\beta}(-\boldsymbol{\Theta}\mathbf{T} \boldsymbol{\Sigma} \mathbf{T}^{*}/4)}{t!},
$$
where, $\re(\beta n/2)> (m-1)\beta/2 - t_{m}$, $\tau = (t_{1}, t_{2},\dots, t_{m})$, $t_{1}\geq t_{2}\geq
\cdots t_{m} \geq 0$, $\sum_{i=1}^{m} t_{i}=t$ and $t_{1}, t_{2},\dots, t_{m}$ are nonnegative integers;
and $\kappa = (k_{1}, k_{2}, \dots, k_{m})$, $k_{1}\geq k_{2}\geq \cdots \geq k_{m} \geq 0$, $k_{1},
k_{2},\dots, k_{m}$ are nonnegative integers.
\end{corollary}
\begin{proof}
Let $\boldsymbol{\Sigma} \in \boldsymbol{\Phi}_{m}^{\beta}$, $\boldsymbol{\Theta} \in
\boldsymbol{\Phi}_{n}^{\beta}$, $\boldsymbol{\mu} \in \mathfrak{L}^{\beta}_{n,m}$, and let $\mathbf{Z}
\sim \mathcal{K}\mathfrak{R}^{\beta, I}_{n \times m} (\kappa,\boldsymbol{0}, \mathbf{I}_{n},
\mathbf{I}_{m})$; and defines $\mathbf{X} = \boldsymbol{\mu}+
\boldsymbol{\Theta}^{1/2}\mathbf{Z}\boldsymbol{\Sigma}^{1/2}$, where $\mathbf{B}^{1/2}$ is the positive
definite square root of $\mathbf{B}$, such that $\mathbf{B}^{1/2}\mathbf{B}^{1/2} = \mathbf{B}$. Then
$\mathbf{X} \sim \mathcal{K}\mathfrak{R}^{\beta, I}_{n \times m} (\kappa,\boldsymbol{\mu},
\boldsymbol{\Theta}, \boldsymbol{\Sigma})$. Therefore
\begin{eqnarray*}
  \phi_{_{\mathbf{X}}}(\mathbf{T}) &=& \E(\etr\{i\mathbf{XT}^{*}\}) = \E(\etr\{i\left(\boldsymbol{\mu}+
    \boldsymbol{\Theta}^{1/2}\mathbf{Z}\boldsymbol{\Sigma}^{1/2})\mathbf{T}^{*}\}\right) \\
  &=& \etr\{i\boldsymbol{\mu}\mathbf{T}^{*}\} \E(\etr\{i
    \boldsymbol{\Theta}^{1/2}\mathbf{Z}\boldsymbol{\Sigma}^{1/2}\mathbf{T}^{*}\}) \\
  &=& \etr\{i\boldsymbol{\mu}\mathbf{T}^{*}\} \E(\etr\{i
    \mathbf{Z}\left(\boldsymbol{\Theta}^{1/2}\mathbf{T}\boldsymbol{\Sigma}^{1/2})^{*}\}\right)\\
  &=& \etr\{i\boldsymbol{\mu}\mathbf{T}^{*}\}
  \phi_{_{\mathbf{Z}}}\left(\boldsymbol{\Theta}^{1/2}\mathbf{T}\boldsymbol{\Sigma}^{1/2}\right).
\end{eqnarray*}
The conclusion is follows as consequence immediately of  Theorem \ref{teocfKR1}.
\end{proof}

\section*{Conclusions}

There is no doubt of the importance of the Riesz distribution from a theoretical point of view and, by
establishing a generalisation of the Wishart distribution in imminent great importance from a practical
point of view. For some time, statisticians have been trying to extend the theory of sampling in
multivariate analysis for samples of a non-normal population. The introduction of Kotz-Riesz distribution
and the distribution of Riesz are an advance in that direction, furthermore, now these results, combined
with those obtained in \citep{dg:15a}, \citep{dg:15b} allow generalise diverse techniques in multivariate
analysis, assuming a Riesz distribution instead an Wishart distribution or equivalently, assuming a
Kotz-Riesz distribution instead a matrix multivariate normal distribution. Finally emphasises that
Kotz-Riesz distribution belongs to the family of left-elliptical distributions, then all the results
obtained for the latter family, as moments, estimation, hypothesis testing, etc., can be quite
\textit{easily} particularised for Kotz-Riesz distribution.

\section*{Acknowledgements}
The author wish to thank the Editor and the anonymous reviewers for their constructive comments on the
preliminary version of this paper. This article was written under the existing research agreement between
the first author and the Universidad Aut\'onoma Agraria Antonio Narro, Saltillo, M\'exico.

\end{document}